\documentclass{article}
\usepackage{amsmath}
\usepackage{amssymb}
\usepackage{epsfig,subfigure,verbatim,epstopdf,graphics}
\usepackage{tabularx}
\usepackage{color}

\newtheorem{theorem}{Theorem}[section]

\newtheorem{corollary}[theorem]{Corollary}

\newenvironment{proof}[1][Proof]{\begin{trivlist}
\item[\hskip \labelsep {\bfseries #1}]}{\end{trivlist}}

\newenvironment{remark}[1][Remark]{\begin{trivlist}
\item[\hskip \labelsep {\bfseries #1}]}{\end{trivlist}}

\newcommand{\qed}{\nobreak \ifvmode \relax \else
      \ifdim\lastskip<1.5em \hskip-\lastskip
      \hskip1.5em plus0em minus0.5em \fi \nobreak
      \vrule height0.75em width0.5em depth0.25em\fi}

\newcommand{\x}{\textbf{x}}
\renewcommand{\v}{\textbf{v}}
\newcommand{\w}{\textbf{w}}
\renewcommand{\k}{\textbf{k}}
\renewcommand{\d}{\partial}
\newcommand{\R}{\mathbb{R}}
\renewcommand{\u}{\textbf{u}}

\newcommand{\eps}{\varepsilon}

\newcommand{\Qt}{\tilde{Q}}

\def \Ghat {\widehat{G}}

\begin{document}

\title{A conservative spectral method for the Boltzmann equation with anisotropic scattering and the grazing collisions limit}

\author{ \and Irene M. Gamba \thanks{Department of Mathematics, The University of Texas at Austin, 2515 Speedway, Stop C1200 Austin, Texas 78712 and ICES, The University of Texas at Austin, 201 E. 24th St., Stop C0200, Austin, TX 78712} \and Jeffrey R. Haack \thanks{Department of Mathematics, The University of Texas at Austin, 2515 Speedway, Stop C1200 Austin, Texas 78712} 
}

\maketitle

\begin{abstract}
 We present the formulation of a conservative spectral method for the Boltzmann collision operator with anisotropic scattering cross-sections. The method is an extension of the conservative spectral method of Gamba and Tharkabhushanam \cite{GamTha09, GamTha10}, which uses the weak form of the collision operator to represent the collisional term as a weighted convolution in Fourier space. The method is tested by computing the collision operator with a suitably cut-off angular cross section and comparing the results with the solution of the Landau equation. We analytically study the convergence rate of the Fourier transformed Boltzmann collision operator in the grazing collisions limit to the Fourier transformed Landau collision operator under the assumption of some regularity and decay conditions of the solution to the Boltzmann equation. Our results show that the angular singularity which corresponds to the Rutherford scattering cross section is the critical singularity for which a grazing collision limit exists for the Boltzmann operator. Additionally, we numerically study the differences between homogeneous solutions of the Boltzmann equation with the Rutherford scattering cross section and an artificial cross section, which give convergence to solutions of the Landau equation at different asymptotic rates. We numerically show the rate of the approximation as well as the consequences for the rate of entropy decay for homogeneous solutions of the Boltzmann equation and Landau equation. 

\end{abstract}



\section{Introduction} \label{sec:intro}

The initial focus of this manuscript was the study of simulating the Boltzmann equation with anisotropic, singular angular scattering cross sections by spectral methods. However, while attempting to verify numerical results based on this method by examining the grazing collision limit of the Boltzmann operator, we found an analytical argument that not only gives an explicit representation of the effect of angular averaging for a family of singular grazing collision angular cross sections, but also the rate of convergence of the grazing collision limit of the the Boltzmann operator to the Landau operator. The bulk of this manuscript will address both the numerical and analytical aspects of this grazing collision limit of the Boltzmann equation in physically relevant regimes, which includes the case of Coulombic intermolecular potential scattering mechanisms.

While numerical methods for solving the Boltzmann equation generally use the assumption of spherical particles with `billiard ball' like collisions, a more physical model is to assume that particles interact via two-body potentials. Under this assumption the Boltzmann equation can be formulated in a very similar manner \cite{Cercignani_Boltz}, but in this case the scattering cross section is highly anisotropic in the angular variable. In many cases, such as the physically relevant case of Coulombic interactions between charged particles, the derivation of the Boltzmann equation breaks down completely due to the singular nature of this scattering cross section. Physical arguments by Landau \cite{Landau37} as well as a later derivation by Rosenbluth et al. \cite{Rosenbluth57} showed that the dynamics of the Boltzmann equation can approximated by a Fokker-Planck type equation when grazing collisions dominate, generally referred to as the Landau or Landau-Fokker-Planck equation. Later work \cite{BobLan75, DesLandau, DegLucLandau,
VillaniThesis, AlexVil} more rigorously justified this asymptotic limit. 

Many numerical methods have been developed for solving the full Landau equation, some stochastic \cite{TakiAbe, Nanbu97} and some deterministic \cite{ParRusTos00}, however very few methods have been developed to compute the Boltzmann equation near this grazing collision limit. The small parameter used to quantify this limit is related to the physical Debye length, which quantifies the distance at which particles are screened from interaction, and a heuristic minimum interaction distance for the grazing collisions assumption to hold. Other non-grazing effects with the Boltzmann equation may remain relevant \cite{FriBook} which makes development of numerical methods based on the Boltzmann equation itself relevant for plasma applications. To our knowledge the only numerical method that makes this distinction explicit is the recently proposed Monte Carlo method for the Landau equation of Bobylev and Potapenko \cite{BobPot}, which grew out of the work of Bobylev and Nanbu \cite{BobNan00}. Pareschi, 
Toscani, and Villani \cite{ParTosVil03} showed that the weights of their spectral Galerkin method for the Boltzmann equation converged to the weights of a similar method for the Landau equation, but neither estimates nor computations were done for the Boltzmann equation near the grazing collisions limit. This work seeks to bridge that gap using the conservative spectral method for the Boltzmann equation developed by Gamba and Tharkabhushanam \cite{GamTha09, GamTha10}.

There are many difficulties associated with numerically solving the Boltzmann equation, most notably the dimensionality of the problem and the conservation of the collision invariants. For physically relevant three dimensional applications the distribution function is seven dimensional and the velocity domain is unbounded. In addition, the collision operator is nonlinear and requires evaluation of a five dimensional integral at each point in phase space. The collision operator also locally conserves mass, momentum, and energy, and any approximation must maintain this property to ensure that macroscopic quantities evolve correctly.

Spectral methods are a deterministic approach that compute the collision operator to high accuracy by exploiting its Fourier structure. These methods grew from the analytical works of Bobylev \cite{Bob88} developed for the Boltzmann equation with Maxwell type potential interactions and integrable angular cross section, where the corresponding Fourier transformed equation has a closed form. Spectral approximations for these type of models where first proposed by Pareschi and Perthame \cite{ParPer96}.
Later Pareschi and Russo \cite{ParRus00} applied this work to variable hard potentials by periodizing the problem and its solution and implementing spectral collocation methods.
 
These methods require $O(N^{2d})$ operations per evaluation of the collision operator, where $N$ is the total number of velocity grid points in each dimension. While convolutions 
can generally be computed in $O(N^d \log N)$ operations, the presence of the convolution weights requires the full $O(N^{2d})$ computation of the convolution,  except for a few special cases such as hard spheres in 3D (and Maxwell molecules in 2D) which can be done with $O(N^{\frac{5}3} \log N)$ in 3D. Spectral methods have advantages over Direct Simulation Monte Carlo Methods (DSMC) in many applications, in particular time dependent problems, low Mach number flows, high mean velocity flows, and flows that significantly deviate from equilibrium. In addition, deterministic methods avoid the statistical fluctuations that are typical of particle based methods. 

Inspired by the work of Ibragimov and Rjasanow \cite{IbrRja02}, Gamba and Tharkabhushanam \cite{GamTha09, GamTha10} observed that the Fourier transformed collision operator takes the form of a weighted convolution and  developed a spectral method based on the weak form of the Boltzmann equation that provides a general framework for computing both elastic and inelastic collisions. Macroscopic conservation is enforced by solving a numerical constrained optimization problem that finds the closest distribution function in $L^2$  in the computational domain to the output of the collision term that conserves the macroscopic quantities. This optimization problem is the approximation of the projection of the Boltzmann solution to the space of collision invariants associated to the corresponding collision operator   \cite{AlonGamTha}. In addition  these methods do not impose periodization on the function but rather assume that solution of the underlying problem on the whole phase space is obtained by the use of the Extension Operator in Sobolev spaces. They are also shown in the space homogeneous, hard potential case, to converge to the Maxwellian distribution with the conserved moments corresponding to the endowed collision invariants. Such convergence and error estimates results heavily rely on the discrete constrained optimization problem (see \cite{AlonGamTha} for a complete proof and details.)

The proposed computational approach is complemented  by the analysis of the approximation from the Boltzmann operator with grazing collisions to the Landau operator by estimating the $L^\infty$-difference  of their  Fourier transforms evaluated on the solution of the corresponding Boltzmann equation, as they both can be easily expressed by  a weighted convolution structure in Fourier space.  We show that this property holds for a family of singular angular scattering cross sections with a suitably cut-off Coulomb potential.  The parameter $0 \le \delta < 2$ corresponds to the strength of the singularity in the angular cross section and the parameter $\eps$ corresponds to the angular cutoff, which gives the \emph{grazing collision limit} as $\eps \to 0$.

The case when the parameter $\delta=0$ corresponds to the classical Rutherford scattering cross section \cite{Ruther11}, and includes an inverse logarithmic term in $\eps$ that ensures the limit.  This value of $\delta$ is critical to obtain the grazing collision limit in the following sense: if $\delta<0$ then the error between the Boltzmann and the Landau operators will not necessarily converge to zero in $\eps$. In addition, for any other value $0<\delta\ < 2$, the rate of convergence of the Boltzmann to the Landau operator is faster in $\eps$. In these sense we can assert that the Rutherford scattering cross section \cite{Ruther11} is the one that contains the weakest possible singularity in the angular cross section for which one can achieve a grazing collision limit to the Landau equation. 
%

These results are shown in  Theorem~\ref{theo1} and are written for the three dimensional case. There we prove that the  $L^\infty$-difference of the Fourier transforms between the Landau operator $Q_L$ and the Boltzmann operator for this family of cross sections  $Q_{b_{\eps}^\delta}$, both  acting on solutions  $f_{\eps}^\delta(v)$  of the Boltzmann equation,
 converges to zero in $\eps$ with rates depending on $\delta$. This requires that the solution $f_{\eps}^{\delta}(v,t)$ satisfies the regularity and decay condition  $\mathcal{F}(f_{\eps, \cdot}^{\delta}\,\tau_{{\u}}f_{\eps}^{\delta})(\zeta, \cdot)< \mathcal{A}(\zeta)/(1+|\u|^{3+a})$ with $a>0$ and
{$\mathcal{A}(\zeta) \le k(1+|\zeta|)^{-3}$,  uniformly in $\R^3$}.
 Our analysis shows the convergence rate in $\eps$ explicitly depends on the parameter $\delta$ that
  quantifies the strength of the non-integrable singularity associated to the collision angular cross section.


In addition, we examine the consequences of this theorem by numerically studying the differences between Rutherford scattering cross section ($\delta=0$), which has logarithmic error in approximating Landau, and the cross section corresponding to $\delta=1$, which better approximates Landau. These numerical results clearly  exhibit the speed of the approximation and decay rate of the entropy functional  for homogeneous solutions of the Boltzmann equation with cut-off Columbic interactions and are benchmarked with  the  solution to the Landau equation, which is independent of the $\eps$ and $\delta$ parameters.



This article is organized as follows. In Section 2, we present the derivation of the spectral formulation of the collision operator for an arbitrary anisotropic scattering cross section. In Section 3 we present the Landau equation and and present also a  broad class of singular angular cross sections formulated by Villani \cite{Villani98} and Bobylev \cite{BobPot} satisfying suitable conditions for the study of grazing collisions limits of the anisotropic Boltzmann equation. We then introduce a family of angular primitives parametrized by $\delta$ to define admissible angular cross sections $b_\eps^\delta(\cos\theta)$ of the scattering angle $\theta$ in order to achieve a grazing collision limit.  This family of angular cross sections includes the screened Rutherford cross section for Coulombic interactions. Then, we prove in  Theorem~\ref{theo1}   the estimates for the difference of the Fourier Transforms of the Landau and Boltzmann grazing collisional operators in terms of the $\eps$ and $\delta$ parameters that yields, as shown in Corollary~\ref{coro1},   
  the rate of asymptotic convergence of solutions of  the Boltzmann collision  to approximate solutions to the  Landau equation for this family cross sections given some condition on the solution of the corresponding Boltzmann equation as described above. In Section 4, we present the details of the numerical method based on this formulation and provide some practical observations on its implementation. In Section 5, we numerically investigate the method's performance for small but finite values of $\eps$, the grazing collision parameter, for the choice of $\delta=0$ and $\delta=1$. We conclude with a discussion of future work in this area.




\section{The space homogeneous Boltzmann equation} \label{sec:BTE}
The space homogeneous elastic Boltzmann equation is given by the initial value problem
\begin{equation}\label{BTE}
\frac{\d}{\d t} f(\v,t) = Q(f,f)(\v,t),
\end{equation}
with
\begin{align*}
\qquad \v \in \R^d, \qquad f(v,0) = f_0(\v) \\
\end{align*}
where $f(\v,t)$ is a probability density distribution in $\v$-space and  $f_0$ is assumed to be at least locally integrable with respect to $\v$. 
 
The collision operator $Q(f,f)$ is a bilinear integral form in $(\v,t)$ given by
\begin{equation}\label{Q_general}
Q(f,f)(\v,t) = \int_{\v_\ast \in \R^d} \int_{\sigma \in S^{d-1}} B(|\v - \v_\ast|,\cos \theta) (f(\v_\ast ')f(\v') - f(\v_\ast)f(\v)) d\sigma d\v_\ast,
\end{equation}
where the velocities $\v', \v_\ast'$ are determined through a given collision rule depending on $\v, \v_\ast$. The positive term of  the integral in \eqref{Q_general} evaluates  $f$ in the pre-collisional velocities that can result in a post-collisional velocity the direction $\v$.  The scattering cross section $B(|\v - \v_\ast|,\cos \theta)$ is a given non-negative function depending on the size of the relative velocity $\u := \v - \v_\ast$ and $\cos \theta = \frac{\u \cdot \sigma}{|\u|}$, where $\sigma$ in the $n-1$ dimensional sphere $S^{n-1}$ is referred to as the scattering direction, which coincides with the direction of the post-collisional elastic relative velocity.

The elastic (or reversible) interaction law written in the scattering direction $\sigma$ is given by
\begin{align}
&\v' = \v + \frac{1}{2}(|\u|\sigma - \u), \qquad \v_\ast ' = \v_\ast - \frac{1}{2}(|\u|\sigma - \u) \label{velocity_interact} \\
& B(|\u|,\cos \theta) = |\u|^\lambda b(\cos \theta)  \, .   \nonumber
\end{align}
The angular cross section function $b(\cos \theta)$ may or may not be integrable with respect to $\theta$; the case when integrability holds is referred to as the Grad cut-off assumption on the angular cross section. 

The parameter $\lambda$ regulates the collision frequency as a function of the relative velocity $|\u|$. This parameter corresponds to the interparticle potentials used in the derivation of the collisional cross section and choices of $\lambda$ are denoted as variable hard potentials (VHP) for $0 < \lambda < 1$, hard spheres (HS) for $\lambda = 1$, Maxwell molecules (MM) for $\lambda = 0$, and variable soft potentials (VSP) for $-3 < \lambda < 0$. The $\lambda = -3$ case corresponds to a Coulombic interaction potential between particles. If $b(\cos \theta)$ is independent of $\theta$ we call the interactions isotropic, e.g., in the case of hard spheres in three dimensions.

\subsection{Spectral formulation for anisotropic angular cross section} \label{sec:spectral_cont}
The key step in our formulation of the spectral numerical method is the use of the weak form of the Boltzmann collision operator \cite{GamTha09}. For a suitably smooth test function $\phi(\v)$ the weak form of the collision integral is given by the double mixing operator
\begin{equation} \label{collision_weakform}
\int_{\R^d} Q(f,f) \phi(\v) d\v = \int_{\R^d \times \R^d \times S^{d-1}} f(\v)f(\v_\ast) B(|\u|,\cos \theta) (\phi(\v') - \phi(\v)) d\sigma d\v_\ast d\v ', 
\end{equation}
If one chooses 
\begin{equation*}
\phi(\v) = e^{-i \zeta \cdot \v} / (\sqrt{2\pi})^d, 
\end{equation*}
then,  \eqref{collision_weakform} is the Fourier transform of the collision integral with Fourier variable $\zeta$:
\begin{align}\label{FourierQ}
\widehat{Q}(\zeta) &= \frac{1}{(\sqrt{2\pi})^d} \int_{\R^d} Q(f,f) e^{-i \zeta \cdot \v} d\v \nonumber \\
&= \int_{\R^d \times \R^d \times S^{d-1}} f(\v)f(\v_\ast) \frac{B(|\u|,\cos \theta)}{(\sqrt{2\pi})^d} (e^{-i \zeta \cdot \v'} - e^{-i \zeta \cdot \v}) d\sigma d\v_\ast d\v \nonumber\\
&= \int_{\R^d} G_b(\u,\zeta) \mathcal{F}[f(\v)f(\v-\u)](\zeta) d\u,
\end{align}
where $\widehat{[\cdot]} = \mathcal{F}(\cdot)$ denotes the Fourier transform and 
\begin{equation} \label{G_eqn}
G_b(\u,\zeta) = |\u|^\lambda \int_{S^{d-1}} b(\cos \theta) \left(e^{-i\frac{\zeta}{2}\cdot( -\u + |\u|\sigma)} - 1\right) d\sigma
\end{equation}
Further simplification can be made by writing the Fourier transform inside the integral as a convolution of Fourier transforms:
\begin{align} \label{Cont_spectral}
\widehat{Q}(\zeta) 
&= \int_{\R^d} \widehat{G_b}(\xi,\zeta) \hat{f}(\zeta - \xi) \hat{f}(\xi) d\xi,
\end{align}
\ 
where the convolution weights $\widehat{G}(\xi,\zeta)$ are given by
\begin{align} \label{Ghat_eqn}
\widehat{G_b}(\xi,\zeta) &= \frac{1}{(\sqrt{2\pi})^d}  \int_{\R^d} G_b(\u,\zeta) e^{-i \xi \cdot u} d\u \\
&=\frac{1}{(\sqrt{2\pi})^d}  \int_{\R^d} |\u|^\lambda e^{-i\xi\cdot\u} \int_{S^{d-1}} b(\cos \theta)  \left(e^{i\frac{\zeta}{2}\cdot( \u - |\u|\sigma)} - 1\right) d\sigma d\u \nonumber
\end{align}
\ 
These convolution weights can be precomputed once to high accuracy and stored for future use. For many collisional models, such as isotropic collisions, the complexity of the integrals in the weight functions can be reduced dramatically through analytical techniques\cite{GamTha09, GamTha10}. However unlike previous work, in this paper we make no assumption on the isotropy of $b$ and derive a more general formula. We remark that this formulation does not separate the gain and loss terms of the collision operator, which is important for obtaining the correct cancellation in the grazing collision limit below.

We begin with $G_b(\u,\zeta)$ and define a spherical coordinate system for $\sigma$ with a pole in the direction of $\u$, i.e. let $\sigma = \cos\theta \frac{\u}{|\u|} + \sin\theta \omega,\, \omega \in S^{d-2}$. We obtain

\begin{align}\label{bessel-w1} 
G_b(\u,\zeta) &=|\u|^\lambda \int_0^\pi \int_{S^{d-2}} b(\cos \theta) \sin\theta  \left(e^{i\frac12 (1-\cos\theta) \zeta \cdot \u}e^{-i\frac12 |\u|\sin\theta (\zeta \cdot \omega)} - 1\right) d\theta d\omega.
\end{align}
\ 
For the remainder of this paper, we will work in three dimensions ($d=3$). We write the unit vector $\sigma$ as $\sigma = \cos\theta \frac{\u}{|\u|} + \sin\theta (\textbf{j} \sin\phi + \textbf{k} \cos\phi)$, where $\textbf{j},\textbf{k}$ are mutually orthogonal vectors with $\u$. Thus the right hand side of \eqref{bessel-w1} becomes
\ 
\begin{align*} 
 |\u|^\lambda \int_0^\pi \int_{\alpha-\pi}^{\alpha + \pi} b(\cos \theta) \sin\theta  \left(e^{i\frac12 (1-\cos\theta) \zeta \cdot \u}e^{-i\frac12 |\u|\sin\theta (\zeta \cdot \textbf{j}\sin\phi + \zeta \cdot \textbf{k} \cos\phi)} - 1\right) d\theta d\phi, 
\end{align*}
for $\alpha$ to be justified below. 

Using the trigonometric identity 
\[ (\zeta \cdot \textbf{j})\sin\phi + (\zeta \cdot \textbf{k})\cos\phi = \sqrt{(\zeta\cdot\textbf{j})^2 + (\zeta\cdot\textbf{k})^2}\sin(\phi + \gamma), \]
for a unique $\gamma \in [-\pi,\pi]$,  the integration with respect to the azimuthal angle  $\phi$ is equivalent to

\begin{align*} 
G_b(\u,\zeta) &=
|\u|^\lambda \int_0^\pi b(\cos \theta)  \sin\theta \left(e^{i\frac12 (1-\cos\theta) \zeta \cdot \u}\int_{\alpha-\gamma-\pi}^{\alpha -\gamma + \pi} e^{-i\frac12 |\u|\sin\theta |\zeta^\perp| \sin\phi}d\phi- 2\pi \right) d\theta, \nonumber \\
&= |\u|^\lambda \int_0^\pi b(\cos \theta)  \sin \theta \left(e^{i\frac12 (1-\cos\theta) \zeta \cdot \u}\int_{\alpha-\gamma-3\pi/2}^{\alpha -\gamma + \pi/2} e^{i\frac12 |\u|\sin\theta |\zeta^\perp| \cos\phi}d\phi- 2\pi \right) d\theta,
\end{align*}
where $\zeta^\perp = \zeta - (\zeta \cdot \textbf{u}/|\u|) \textbf{u}/|\u|$.  Finally, let $\alpha = \gamma + \pi/2$, then by symmetry we obtain

\begin{align} 
G_b(\u,\zeta) &= |\u|^\lambda \int_0^\pi b(\cos \theta) \sin\theta \left(e^{i\frac12 (1-\cos\theta) \zeta \cdot \u} 2\int_{0}^{\pi} e^{i\frac12 |\u|\sin\theta |\zeta^\perp| \cos\phi}d\phi- 2\pi \right) d\theta \nonumber\\
&= 2\pi |\u|^\lambda \int_0^\pi b(\cos \theta) \sin\theta \left(e^{i\frac12 (1-\cos\theta) \zeta \cdot \u} J_0\left(\frac{|\u|\sin\theta|\zeta^\perp|}{2}\right) - 1 \right) d\theta, \label{SphG}
\end{align}
where $J_0$ is the Bessel function of the first kind (see \cite{AbrSte} 9.2.21). Note that for the isotropic case the angular function  $b(\cos \theta)$ is constant and thus $\zeta$ can be used instead of $\u$ as the polar direction for $\sigma$, resulting in an explicit expression involving a sinc function \cite{GamTha09}.

Next, we take  $\widehat{G_b}$ to be the Fourier transform of $G_b$. Note that by symmetry $\widehat{G_b}$ is real valued, thus this transform is taken on a ball centered at 0 in order to ensure that this symmetry is maintained while computing them numerically.

Then, the convolution weights $\Ghat_b(\zeta,\xi)$ from \eqref{Ghat_eqn}, written in 3 dimensions, are computed as follows
\begin{align*}
\Ghat_b(\xi,\zeta) &= 2\pi \int_{\R^3} |\u|^\lambda e^{-i\xi \cdot \u} \int_0^\pi b(\cos\theta)\sin\theta  \\
\times & \left[e^{\frac{i\zeta}{2} \cdot \u(1-\cos\theta)} J_0\left(\frac12 |\u||\zeta^\perp|\sin\theta\right) - 1\right] d\theta d\u \\
&=2\pi \int_0^\infty  \int_{S^2} r^{\lambda+2} \int_0^\pi b(\cos\theta)\sin\theta \\
\times & \left[e^{-ir(\xi - \frac{\zeta}{2}(1-\cos\theta) )\cdot \eta}J_0\left(\frac12 r|\zeta^\perp|\sin\theta\right)  -  e^{-ir\xi \cdot \eta} \right] d\theta d\eta dr. \label{Ghat_Aniso}
\end{align*}

We now take $\zeta$ to be the polar direction for the spherical integration of $\eta$ and use that $\widehat{G}_b$ is real-valued to obtain

\begin{align}
\Ghat_b(\xi,\zeta) &=
4\pi^2 \int_0^\infty  r^{\lambda+2} \int_0^\pi \int_0^\pi b(\cos\theta)\sin\theta\sin\gamma J_0\left(r\left|\xi - \frac{\xi\cdot\zeta}{|\zeta|^2}\zeta\right|\sin\gamma\right)  \nonumber \\
 & \times \Bigg[\cos\left(r(\xi - \frac{\zeta}{2}(1-\cos\theta) )\cdot \frac{\zeta}{|\zeta|}\cos\gamma\right) J_0\left(\frac12 r|\zeta|\sin\gamma\sin\theta\right) \nonumber \\
 &- \cos\left(r\xi \cdot \frac{\zeta}{|\zeta|}\cos\gamma\right) \Bigg] d\theta d\gamma dr, 
\end{align}
where $\gamma$ is the polar angle for the $\eta$ integration.


\section{The grazing collisions limit and convergence to the Landau collision operator}

\subsection{The Landau collision operator}
The Landau collision operator describes binary collisions that only result in very small deflections of particle trajectories, as is the case for Coulomb potentials between charged particles \cite{Ruther11}. This can be shown to be an approximation of the Boltzmann collision operator in the case where the dominant collision mechanism is that of grazing collisions. The operator is given by
\begin{equation}\label{FPL}
Q_L (f,f) = \nabla_\v \cdot \left( \int_{\mathbb{R}^3} |\u|^{\lambda+2} (I - \frac{\u\otimes \u}{|\u|^2}) (f(\v_\ast)\nabla_\v f(\v) - f(\v) (\nabla_\v f)(\v_\ast)) d \v_\ast\right) \,,
\end{equation}
and the weak form of this operator is given by \cite{ParTosVil03}
\begin{align*} 
\int_{\R^3} Q_L(f,f) \phi(\v) dv& = \int_{\R^3} \int_{\R^3} f(\v) f(\v_\ast) \nonumber \\
\times & \left(-4|\u|^\lambda \u \cdot \nabla \phi + |\u|^{\lambda+2}\left(I - \frac{\u \otimes \u}{|\u|^2}\right) : \nabla^2 \phi\right) d\v d\v_\ast\, ,
\end{align*}
where $\nabla^2 \phi$ denotes the Hessian of $\phi$ and `$:$' is the matrix double dot product.

As done in the Boltzmann case above, we choose $\phi$ to be the Fourier basis functions and obtain after some calculation
\begin{align}\label{FPL-Fourier}
\widehat{Q_L}(\zeta) &= \frac{1}{(2\pi)^{3/2}} \int_{\R^3}  \mathcal{F}\{f(\v)f(\v-\u)\}(\zeta) \Bigg(4 i |\u|^\lambda (\u \cdot \zeta)- |\u|^{\lambda+2} |\zeta^\perp|^2 \Bigg) d\u,\end{align}
where $\zeta^{\perp} = \zeta - (\zeta \cdot \u)/ |\u|^2\ \u$, the orthogonal component of $\zeta$ to $\u$. Thus the weight function $G_L(\u,\zeta)$ in terms of \eqref{G_eqn} is now given by
\begin{equation} 
G_L(\u,\zeta) = |\u|^\lambda (4i(\u\cdot\zeta) - |\u|^2|\zeta^\perp|^2) \, . \label{GFPL}
\end{equation}
The $\widehat{G_L}$ used in the final computation is the Fourier transform of $G_L$ with respect to $\u$, but we will work with this representation to make the convergence analysis below more clear.

\subsection{The grazing collisions limit}

To show that the spectral representation of Boltzmann operator is consistent with this form of the Landau operator, we must take the grazing collisions limit within this framework. 
To obtain this limit, it is enough to assume that the cross section satisfies the following properties.

Let $\eps > 0$ be the small parameter associated with the grazing collisions limit and $\delta$ be a parameter associated with the strength of the singularity at $\theta = 0$, to be explained below. A family of cross sections $b^\delta_\eps(\cos\theta)$ represents grazing collisions if \cite{AlexVil, BobPot}:
%
%
\begin{align}\label{ruth-like}
\bullet\qquad  &\lim_{\eps \to 0} 2\pi \int_0^{\pi} b_\eps^\delta(\cos \theta) \sin^2(\theta/2)\,
\sin\theta d \theta =\Lambda_0 <\infty, \qquad \Lambda_0 > 0 \nonumber\\
\ &\ \nonumber\\
\bullet\qquad & 2\pi \int_0^{\pi} b_{\eps}^\delta(\cos \theta) (\sin(\theta/2))^{2+k}\, 
\sin\theta d \theta {\to_{_{\eps\to 0}}} 0  \qquad \mathrm{for}\ \ k > 0 \, . \\
\ &\ \nonumber\\
\bullet\qquad &\forall \theta_0 > 0, \ b_\eps^\delta(\cos(\theta)) {\to_{_{\eps\to 0}}} 0 \text{ uniformly on }\ \theta > \theta_0. \nonumber
\end{align}

These conditions are sufficient show that the  collisional integral operator converges to the Landau operator at a rate that depends on the choice of the angular function  $b_\eps^\delta(\cos \theta)$, independently of $\eps$ and $\delta$, provided the solution $f_\eps^\delta$ of the Boltzmann equation for grazing collisions (1,2,3) with 
\eqref{ruth-like}  satisfies some regularity and decay at infinity,   as it will be shown in Theorem~\ref{theo1}.

The most significant and perhaps physically meaningful example family of cross sections that satisfy these conditions can be generated from Rutherford scattering, corresponding to a family $b_{\eps}(\cos\theta)$ given by
\begin{equation}  
b_\eps(\cos \theta) \sin\theta =   \frac{ \sin\theta }{-\pi \log(\sin(\eps/2))\sin^4(\theta/2)} 1_{\theta \ge \eps} .\label{beps_ruth}
\end{equation}

\begin{remark}[Remark 1:]\label{rem1} We note that the logarithmic term that appears here is the \emph{Coulomb logarithm} originally derived by Landau \cite{Landau37}, where $\eps$ is proportional to the ratio between the mean kinetic energy of the gas and the Debye length. As will be observed later, this rescaling of the cross section is required in order to take the limit $\eps \to 0$, as the form of the Landau equation we are using \eqref{FPL} does not have the Coulomb logarithm. Without this rescaling the leading order term of the collision operator would be the Landau equation \eqref{FPL} simply multiplied by $\log(\sin(\eps/2))$, and the remainder terms would also be multiplied by this factor. 
\end{remark}

Another angular cross section that satisfies  conditions \eqref{ruth-like} is given by 
\begin{equation}
b_\eps(\cos \theta) \sin\theta =   \frac{8 \eps }{\pi \theta^4} 1_{\theta \ge \eps}     \label{beps} ,
\end{equation}
which we will refer to as the $\eps$-linear cross section. While this cross section is not physically motivated, it is useful for numerical convergence studies. Other angular cross sections that satisfy  conditions \eqref{ruth-like}  have been used in DSMC methods for computing the Landau equation; for an overview see \cite{BobPot}.

In fact it is possible to identify a large family of possible angular function choices corresponding to two body interaction potentials that includes both the Coulombic case \eqref{beps_ruth} and the $\eps$-linear one \eqref{beps}, the former being the critical case for the grazing collision limit. 

\medskip

\subsection{A family of angular cross sections for long range interactions}

We next introduce a more general way to define a family of angular cross sections that will satisfy conditions  \eqref{ruth-like}. For this purpose we introduce the functions $H(x)$ and $C(x)$  as the primitives of
%
%
 \begin{align}\label{HC-funct}
 H'(x)= b(1-2x^2) x^3 \qquad\qquad\textrm{and}\qquad\qquad C'(x)=b(1-2x^2) x^5 \, .
 \end{align}
for a given (non cut off) cross section $b(\cos \theta)$.

These are related to the grazing conditions by setting $x = \sin(\theta/2)$. Indeed, for $H(x)$ we have 
\begin{align}\label{H-funct}
 H'(x) dx &= b(1-2x^2) x^3 dx = \frac12 b(\cos\theta) \sin^3(\theta/2) \cos(\theta/2) d\theta \nonumber\\
 &=  \frac14 b(\cos\theta) \sin^2(\theta/2) \sin(\theta)  d\theta = \frac18 b(\cos\theta) (1-\cos\theta)\sin(\theta)  d\theta \, .
\end{align}
Note that using this $H$ function we have
\[ \int_\eps^\pi b(\cos\theta) \sin^2(\theta/2) \sin \theta d\theta = 4 \int_{\sin(\eps/2)}^1 H'(x) dx = 4(H(1) - H(\sin(\eps/2))) \]

We similarly define
\begin{align}\label{C-funct}
 C'(x) dx &= b(1-2x^2) x^5 dx = \frac12 b(\cos\theta) \sin^5(\theta/2) \cos(\theta/2) d\theta \nonumber\\
 &= \frac 14 b(\cos\theta) \sin^4(\theta/2) \sin(\theta) d\theta= \frac1{16} b(\cos\theta) \sin(\theta) (1-\cos\theta)^2 d\theta\, ,
\end{align}
for convenience, as it will arise in the proof of the grazing limit.

In order to satisfy the conditions of the second and third bullets for the grazing limit \eqref{ruth-like}, it is sufficient that the angular function $b(\cos\theta)$ is singular enough such that 
\begin{align}\label{HC-conditions}
 &\lim_{\eps \to 0} \frac{1}{H(\sin(\eps/2))} =0  \qquad
 \text{and  }\\
& \left\{|H(1)|, |C(1)|, \sup_{\eps > 0} |C(\sin(\eps/2))| \right\} \le \Gamma. \nonumber
\end{align}
for some constant $\Gamma$ depending only on $b$.

Using these just introduced definitions, the $\eps$-dependent angular cross section with a short range cut-off can be written in terms of the $H$ function from \eqref{HC-funct} as follows 
\begin{align}\label{H-cross}
 b_\eps(\cos\theta) \sin\theta d\theta &= -\frac1{2\pi H(\sin(\eps/2))} b(\cos\theta) \sin \theta \, 1_{\theta\geq \eps}\,  d\theta \nonumber\\
 &= -\frac4{2\pi H(\sin(\eps/2))} \frac{H'(x)}{x^2} 1_{x \geq \sin(\eps/2)}\, dx.
\end{align}
Note that by construction, this cross section clearly satisfies the third grazing limit condition in \eqref{ruth-like}. It also satisfies the first grazing limit condition:
\begin{align*} 2\pi \int_0^\pi b_\eps(\cos\theta) \sin^2(\theta/2) \sin\theta d\theta &= -\frac{1}{H(\sin(\eps/2))} \int_\eps^\pi b(\cos\theta) \sin^2(\theta/2) \sin\theta d\theta \nonumber \\
&=  -\frac{4H(1)}{H(\sin(\eps/2))}  + 4. \\
\end{align*}
For the second grazing condition, note that 
\[ \int_\eps^\pi b_\eps(\cos\theta) \sin^{2+k} (\theta/2) \sin\theta d\theta = -\frac{4}{2\pi H(\sin(\eps/2))} \int^1_{\sin(\eps/2)} x^k H'(x) dx, \qquad k > 0 \]
thus any result will be less singular than the $H(\sin(\eps/2))$ as $\eps \to 0$.

\bigskip

\noindent {\bf A $\delta$-family of admissible angular singularities: }  One can see that when the angular function $b(\hat{u} \cdot \sigma)$ takes the form
\begin{equation}\label{H-cross-2}
 b(\cos \theta) := b^\delta(\cos \theta)=\frac1{\sin^{4+\delta}(\theta/2)}\, ,
\end{equation}
then, after  introducing the $\eps$ and $\delta$ reference parameters in the notation of the  \emph{$\eps$-grazing and $\delta$-singular} angular function $b_\eps^\delta$,  the function $H_\delta(x)$  can be explicitly computed from the area differential of the angular part of the differential cross section \eqref{H-cross}
\begin{align}\label{H-cross-3}
 b_\eps^\delta(\hat{u} \cdot \sigma)d\sigma &= -\frac1{2\pi H_\delta(\sin(\eps/2))} b^\delta(\cos\theta) \, \sin(\theta)\,  1_{\theta\geq \eps}\, d\theta d\omega  \nonumber \\
 &= -\frac1{2\pi H_\delta(\sin(\eps/2))}  \frac1{\sin^{4+\delta}(\theta/2)} \sin(\theta) 1_{\theta\geq \eps}\, d\theta d\omega\nonumber \\
&= -\frac4{2\pi H_\delta(\sin(\eps/2))}   \frac1{x^{1+\delta}} \frac{1}{x^2} \, 1_{x\geq \sin(\eps/2)} dx d\omega .
\end{align}
{Thus, also equating the last term above to the right hand side of relation \eqref{H-cross}, one can explicitly calculate 
$H_\delta(x)$ as the antiderivative of $x^{-(1+\delta)}$, so it has the form
\begin{align}\label{H-cross-4}
H_\delta(x) = -\frac{x^{-\delta}}{\delta},\ \ \text{for} \ \delta>0 \qquad \text{and} \qquad H_0(x)= \log x \, ,\ \ \text{for}\  \delta=0\, .
\end{align}}
Similarly, the corresponding function $C_\delta$, as defined in \eqref{C-funct}, satisfies
\begin{align}\label{C-funct-2}
 C_\delta'(x) dx &=  \frac1{\sin^{4+\delta}(\theta/2)} \sin^{5}(\theta/2) \cos(\theta/2) d\theta \nonumber\\
 & = 2x^{1-\delta} dx \, .
\end{align}

The choice of the exponent $\delta$ must be done in order to satisfy the third bullet condition \eqref{ruth-like}, i.e.
conditions \eqref{HC-conditions} for both $H_\delta(x)$ and $C_\delta(x)$.

\medskip

The case $\delta=0$ yields the Rutherford cross section where
\begin{align}\label{HC-funct-2}
 H_0(x)= \log x \ \ \  \ \text{and} \,\,\,\, C_0(x) = x^2  \,.
\end{align}
These functions satisfy conditions \eqref{HC-conditions}, as 
\begin{equation*}
 \lim_{\eps\to 0} -\frac1{H_0(\sin(\eps/2))}=
\lim_{\eps\to 0} -\frac1{\log(\sin(\eps/2))}= 0.
\end{equation*}

For the $\delta > 0$ case, the $H_\delta(x)$ and $C_\delta(x)$ functions become
\begin{align}\label{HC-funct-3}
 H_\delta(x)= - \frac{x^{-\delta}}{\delta} \ \ \  \ \text{and}  \ \ \ \ \ C_\delta(x) = \frac{2x^{2-\delta}}{2-\delta} .
\end{align}
These two functions also satisfy conditions \eqref{HC-conditions}, as   
\begin{equation*}
 \lim_{\eps\to 0} -\frac1{ H_\delta(\sin(\eps/2))}=
\lim_{\eps\to 0} \frac1 {\sin^{-\delta}(\eps/2)}= \lim_{\eps\to 0}\,  \frac1{2^\delta}\, \eps^{\delta}=0.
\end{equation*}
Finally, notice that the case $\delta=1$ corresponds to the  the $\eps$-linear cross section \eqref{beps}, as $\sin(\theta/2) \approx
\theta/2$ as $\theta \to 0$.  
In this case we have $C_1(x) =  2x$, and thus $C_1(\sin(\eps/2)) =  2\sin(\eps/2)$, satisfying \eqref{HC-conditions}.

The critical case of $\delta = 0$  corresponds to the Rutherford scattering \eqref{beps_ruth}, for which the Landau limit would be possible. 
Clearly,  this case  is the smallest value of the exponent in the singularity
of  the cross section written in negative powers of $\sin(\theta/2)$  such that  the bullet conditions \eqref{ruth-like} for the grazing collision limit are satisfied. 
In this sense \emph{the Coulombic potential case \eqref{beps_ruth} 
 is the critical case for which  the Boltzmann operator can converge to the Landau operator}.

In addition, this approach breaks down when $\delta\ge 2$ as condition  \eqref{HC-conditions} would not be satisfied on $C_\delta(\sin(\eps/2))$.  This value of $\delta$ is the critical one at which more terms in the Taylor expansion for the angular cross-section contain singularities, and the next term of expansion of $G_b(\u,\zeta)$ would need a similar treatment for $C_\delta(x)$ as was done for $H_\delta(x)$ (see the first terms of the expansions in equations \eqref{G13} and \eqref{G16}.)

\bigskip
\subsection{The grazing collision approximation Theorem in three dimensions}

In the following  theorem  we estimate the difference of the grazing collision limit 
for the Boltzmann solutions evaluated at the collisional integral and Landau operators for a class of cross sections given by the general form of  \emph{$\eps$-grazing and $\delta$-singular} angular cross sections satisfying \eqref{H-cross} and \eqref{H-cross-3}.

\medskip

We begin by taking a look at the grazing collisions limit for angular cross sections satisfying conditions  \eqref{H-cross}, and all related conditions for the functions $H$ and $C$ as defined in the previous section.
%
%
%

\begin{theorem}\label{theo1}
 Assume that $f_\eps^\delta$ satisfies
\begin{equation} \label{fAssump}
|\mathcal{F}\{ f_\eps^\delta (\v, t) f_\eps^\delta(\v-\u)\}(\zeta)| \le  \frac{A(\zeta, t)}{1 + |\u|^{3+a}}, 
\end{equation}
with $A(\zeta, t)$ uniformly bounded by {$k(1+|\zeta|)^{-3}$}, $k$ constant, and $a > 0$. We also assume that the angular scattering cross section  $b(\cos \theta)$ satisfies  conditions in \eqref{H-cross} related to the $H_\delta$ function in \eqref{H-cross-4} satisfying conditions \eqref{H-funct} and \eqref{HC-conditions}, with $0 \le \delta < 2$
and $\lambda=-3$.

Then the rate of convergence of the Boltzmann collision operator with grazing collisions  to the Landau collision operator is given by
\begin{align}\label{FTAssump}
 &\|\widehat{Q_L}[f_\eps^\delta] - \widehat{Q_{b_{\eps}^\delta}}[f_\eps^\delta] \|_{L^{\infty}} \  \le  O\left(\frac{  \left | 1 +  (| \log(\sin(\eps/2))|  - 1)  \, 1_{\{\delta=1\}} \right |} { \left | H_\delta(\sin(\eps/2)) \right |}
\right)\ \to_{\eps\to 0 } 0 \, .
\end{align}

\end{theorem}

\bigskip

From this Theorem, the following corollary follows easily by using the assumption that $f_\eps^\delta$ solves the Boltzmann equation for the $\eps,\delta$ family of admissible grazing collision cross sections. Indeed setting $\frac{\d}{\d t} f_\eps^\delta={Q_{b_{\eps}^\delta}}[f_\eps^\delta]$, taking their Fourier transforms and replacing into \eqref{FTAssump} one obtains the following estimate for approximate solutions to the Landau equation.

\begin{corollary}\label{coro1}
Under the conditions of Theorem~\ref{theo1} the following approximation holds
\begin{align}\label{coroll}
 &\| \frac{\partial}{\partial t} \widehat{f_\eps^\delta} - \widehat{Q_L}[f_\eps^\delta]  \|_{L^{\infty}} \  \le  O\left(\frac{  \left | 1 +  (| \log(\sin(\eps/2))|  - 1)  \, 1_{\{\delta=1\}} \right |} { \left | H_\delta(\sin(\eps/2)) \right |}
\right)\ \to_{\eps\to 0 } 0 \, .
\end{align}
\end{corollary}

\begin{remark}[Remark 2:]\label{rem2} Assumption  \eqref{fAssump} means that $f_\eps^\delta$ has at least third order derivatives in $v$ as well as strong decay in $v$, and in fact  rapidly decreasing functions in the Schwarz class satisfy this assumptions. 

 In addition, it is easy to see that  the assumption \eqref{fAssump} is well satified for  $f(t,\v)$ being a Maxwellian distribution in $\v$ space. Indeed take, for ease of presentation,  $f = e^{-|\v|^2/2}$. Then 
\begin{align}
\mathcal{F} \{f(\v) f(\v-\u) \}(\zeta) &= \int_{\R^3} e^{-|\v|^2/2} e^{-|\v-\u|^2/2} e^{-i \zeta\cdot \v} d\v \nonumber \\
&= e^{-|\u|^2/2} \int_{\R^3} e^{-(\v\cdot\v - \v\cdot\u)} e^{-i\zeta\cdot\v} d\v \nonumber \\
&= e^{-|\u|^2/4} \int_{\R^3} e^{-|\v - \u/2|^2} e^{-i\zeta\cdot\v} d\v \nonumber \\
&= e^{-3|\u|^2/4} \int_{\R^3} e^{-|\w|^2} e^{-i\zeta\cdot\w} d\w \nonumber \\
&= \frac{1}{\sqrt{2}}e^{-3|\u|^2/4} e^{-|\zeta|^2/4} \le \frac{{(1+|\zeta|)^{-3}}}{\sqrt{2}(1 + |\u|^{3+a})}
\end{align}
\end{remark}

\begin{remark}[Remark 3:]\label{rem3}  We observe that, as expected from the result of Theorem~\ref{theo1},  the decay rate to equilibrium for the Rutherford $\eps$-logarithmic cross section  \eqref{beps_ruth} is much faster than the one for the  $\eps$-linear cross section  \eqref{beps}, and the latter one actually mimics the entropy decay rate of the Landau equation. 
This fact is well illustrated in Section 5 where we show the numerically computed entropy decay associated to the solution of the initial value problem for Boltzmann with Rutherford cross section \eqref{beps_ruth}. This is in fact an expected observation, as we explain below.
\end{remark}

\bigskip

\begin{proof}{\bf of Theorem~\ref{theo1}.}
With this angular cross section and $\lambda = -3$, the calculation for the weight function $G_{b_\eps}(\zeta,\u)$ can be computed by Taylor expanding the exponential term in \eqref{G_eqn} to obtain:
\begin{align}\label{G_beps}
G_{b_\eps^\delta}(\zeta,\u) &= |\u|^{-3} \int_{S^2}  b_\eps^\delta(\hat{u} \cdot \sigma)  (e^{-i\frac{\zeta}{2} \cdot (|\u|\sigma - \u)} - 1) d\sigma \nonumber \\
                                  &= |\u|^{-3} \int_{S^2} b_\eps^\delta(\hat{u} \cdot \sigma)  \Bigg(e^{i\left(\frac{\u \cdot \zeta}{2} - |\u| \frac{\zeta \cdot \sigma}{2}\right)} - 1\Bigg) d\sigma \,  \\
                                  &= |\u|^{-3} \int_{S^2} b_\eps^\delta(\hat{u} \cdot \sigma) \Bigg[ i \left(\frac{\u \cdot \zeta}{2} - |\u| \frac{\zeta \cdot \sigma}{2}\right) \nonumber \\
                                  & - \frac12 \left(\frac{\u \cdot \zeta}{2} - |\u| \frac{\zeta \cdot \sigma}{2}\right)^2  - i e^{ic} \left(\frac{\u \cdot \zeta}{2} - |\u| \frac{\zeta \cdot \sigma}{2}\right)^3 \Bigg] d\sigma\, , \nonumber
\end{align}
for some $c$ such that $0 < |c| < \left|\frac{\u \cdot \zeta}{2} - |\u| \frac{\zeta \cdot \sigma}{2}\right|$.

We then define $\sigma$ in terms of the polar and azimuthal angles $\theta$ and  $\phi$, respectively associated to the change of coordinates $\sigma = \cos\theta \frac{\u}{|\u|} + \sin\theta \omega$, with now $\omega \in S^{1}$ used in Section 2.1:
As a consequence,   the following representation for the weight function  $G_{b_\eps^\delta}(\zeta,\u)$ holds
\begin{align}\label{G_beps2}
&G_{b_\eps^\delta}(\zeta,\u)  = |\u|^{-3} \int_0^{\pi} \int_0^{2\pi} b_\eps^\delta(\cos \theta) \sin\theta \Bigg[ i \left(\frac{(\u \cdot \zeta)(1-\cos\theta)}{2} - |\u|\frac{\zeta \cdot \hat{\zeta}^\perp}{2} \sin\theta \sin\phi \right) \nonumber \\
                                  &\qquad\qquad- \frac12 \left(\frac{(\u \cdot \zeta)(1-\cos\theta)}{2} - |\u|\frac{\zeta \cdot \hat{\zeta}^\perp}{2} \sin\theta \sin\phi \right)^2 \nonumber \\
                                  &\qquad\qquad- \frac{i e^{ic}}{6} \left(\frac{(\u \cdot \zeta)(1-\cos\theta)}{2} - |\u|\frac{\zeta \cdot \hat{\zeta}^\perp}{2} \sin\theta \sin\phi \right)^3 \Bigg] d\phi d\theta \nonumber \\
                                  &\ \ = |\u|^{-3} \int_0^{\pi} \int_0^{2\pi} b_\eps^\delta(\cos \theta) \sin\theta \Bigg[ i \left( (\u \cdot \zeta)\sin^2(\theta/2) - |\u|\zeta \cdot \hat{\zeta}^\perp \sin(\theta/2) \cos(\theta/2) \sin\phi  \right) \nonumber \\
                                  &\qquad - \frac12 \left( (\u \cdot \zeta)\sin^2(\theta/2) - |\u|\zeta \cdot \hat{\zeta}^\perp \sin(\theta/2) \cos(\theta/2) \sin\phi  \right)^2  \nonumber \\
                                  &\qquad - \frac{i e^{ic}}{6} \left( (\u \cdot \zeta)\sin^2(\theta/2) - |\u|\zeta \cdot \hat{\zeta}^\perp \sin(\theta/2) \cos(\theta/2) \sin\phi  \right)^3 \Bigg] d\phi d\theta \nonumber \\
                                  &:= G_{b_\eps^\delta,1} + G_{b_\eps^\delta,2} + G_{b_\eps^\delta, 3}\, , 
\end{align}
where $\hat{\zeta}^\perp$ is the unit vector in the direction of the part of $\zeta$ that is orthogonal to the pole $\u$, which arises from this choice of spherical coordinates, and note that $\zeta \cdot \hat{\zeta}^\perp = |\zeta^\perp|$. We stress that this expansion occurs in the convolution weights in this formulation rather than the distribution function as is done in other derivations of the grazing collisions limit.

Next,  we use the form of the angular  cross section function $b=b_\eps^\delta(\cos\theta)$  as defined through \eqref{H-cross},  \eqref{H-cross-2} and \eqref{H-cross-4};
 and  examine the result arising from the first two terms of the expansion to obtain
\begin{align}\label{G12}
& G_{b_{\eps}^\delta,1} + G_{b_{\eps}^\delta,2}=\nonumber \\
&\ \ \ \frac{|\u|^{-3}}{-2\pi H_\delta(\sin(\eps/2))} \int_\eps^\pi \int_0^{2\pi}  b_{\eps}^\delta(\cos\theta)  \sin(\theta) \Bigg( i\sin^2(\theta/2) (\u \cdot \zeta) \nonumber \\
& \ \ \ \ \ \ \ - i|\u||\zeta^\perp| \sin(\theta/2)\cos(\theta/2) \sin\phi - \frac{\sin^2(\theta/2)}{2} \nonumber \\
&\ \ \ \ \ \ \ \ \ \ \times \Big((\u\cdot\zeta)^2\sin^2(\theta/2) - 2|\u||\zeta^\perp| (\u\cdot\zeta)\sin(\theta/2)\cos(\theta/2) \sin\phi \nonumber \\
&\ \ \ \ \ \ \ \ \ \ \ \ \  + |\u|^2|\zeta^\perp|^2 \cos^2(\theta/2)\sin^2 \phi \Big)\Bigg) d\phi d\theta.
 \end{align}
 Then integrating in $\phi$ and recalling the definitions of $H_\delta(x)$ and $C_\delta(x)$ from \eqref{HC-funct} to obtain
\begin{align}\label{G13}
&\ =\frac{|\u|^{-3}}{- 2H_\delta(\sin(\eps/2))} \int_\eps^\pi  b(\cos\theta)  \sin(\theta) \Bigg( 2 i\sin^2(\theta/2) (\u \cdot \zeta) \\
&- (\u\cdot\zeta)^2\sin^4(\theta/2) - \frac12 |\u|^2|\zeta^\perp|^2\cos^2(\theta/2)\sin^2(\theta/2) \Bigg) d\theta \nonumber \\
 &\ =\frac{|\u|^{-3}}{-H_\delta(\sin(\eps/2))}  
 \int_{\sin(\eps/2)}^1 \Bigg[ \Big(4 i(\u \cdot \zeta) H'_\delta(x)  - 2(\u\cdot\zeta)^2  \, C'_\delta(x) \Big)  
  - |\u|^2  |\zeta^\perp|^2  \Big(H'_\delta(x) - C'_\delta(x) \Big)\Bigg]  dx  \nonumber \\
 &\ = \frac{|\u|^{-3}}{-H_\delta(\sin(\eps/2))}  \Big[ 4 i(\u\cdot\zeta)\Big( H_\delta(1) - H_\delta(\sin(\eps/2)) \Big) - 
  2 (u\cdot\zeta)^2  \Big( C_\delta(1) - C_\delta(\sin(\eps/2)) \Big) \nonumber\\
& \ \ \ \ - |\u|^2|\zeta^\perp|^2 (H_\delta(1) - C_\delta(1) 
 - H_\delta(\sin(\eps/2)) + C_\delta(\sin(\eps/2))\, )\Big] \nonumber 
 \end{align}
 
\noindent We now invoke the properties of functions $H_\delta(x)$ and $C_\delta(x)$ defined in \eqref{H-funct} and \eqref{C-funct}, respectively,  
 satisfying conditions \eqref{HC-conditions} as well as the identities. Thus, replacing 
 in the last terms of the previous identity in \eqref{G13}, yields
 \ 
 \begin{align}\label{G16} C_\delta(\sin(\eps/2))
 &G_{b_{\eps}^\delta,1} + G_{b_{\eps}^\delta,2}
 = |\u|^{-3}\left(4i(\u\cdot\zeta) - |\u|^2|\zeta^\perp|^2\right) \\
&\ \ - \frac{|\u|^{-3}}{H_\delta(\sin(\eps/2)) }\, \left(  4 i(\u\cdot\zeta) H_\delta(1) - 2 (\u\cdot\zeta)^2 (C_\delta(1)- C_\delta(\sin(\eps/2)))\right. \nonumber\\
&\ \ \ \ \  \ \ \ \ \left. - |\zeta^\perp|^2|\u|^2 {(H_\delta(1) - C_\delta(1) +  C_\delta(\sin(\eps/2)) )}\right) \nonumber 
 \end{align}

Note that the first term is exactly the weight derived for the Landau operator above \eqref{GFPL}.  The remaining terms, 
having the coefficients $C_\delta(1), H_\delta(1)$ and $ C_\delta(\sin({\eps}/2))$ bounded, will vanish as $\eps \to 0$. 
Indeed, defining
\begin{align}\label{Gtilde}
\tilde{G}(\zeta,\u):=    G_{b_{\eps}^\delta}(\zeta,\u)  - G_L(\zeta,\u).
\end{align}
as the deviation from the Landau weights, one obtains
\begin{align}\label{LandauConverge}
 \widehat{Q_{b_{\eps}}}[f_\eps^\delta](\zeta)   &= \widehat{Q_L}[f_\eps^\delta](\zeta) + 
 \int_{\R^3} \mathcal{F}\{f_\eps^\delta(v)f_\eps^\delta(v-u)\}(\zeta)  \tilde{G}(\zeta,\u) du \nonumber\\
 &\ \ = \widehat{Q_L}[f_\eps^\delta](\zeta) +  \int_{\R^3} \mathcal{F}\{f_\eps^\delta(v)f_\eps^\delta(v-u)\}(\zeta) \\
& \ \ \ \ \times \Bigg[\frac{|\u|^{-3}}{H_\delta(\sin(\eps/2)) }\, \left(  2 i(\u\cdot\zeta) H_\delta(1) - (\u\cdot\zeta)^2 (C_\delta(1)- C_\delta(\sin(\eps/2)))\right. \nonumber\\
&\ \ \ \ \ \ \ \ \ \ \  \left.-|\zeta^\perp|^2|\u|^2 {(H_\delta(1) - C_\delta(1) +  C_\delta(\sin(\eps/2)) )}\right)  +  G_{b_{\eps}^\delta,3}(\zeta,\u) \Bigg] du. \  \nonumber
\end{align}

Thus, we need to control $\tilde{G}(\zeta,\u)$ in order to estimate the convergence rate in \eqref{LandauConverge}.
We notice that the leftover terms from $G_{b_\eps^\delta,1},G_{b_\eps^\delta,2}$ are controlled by
\begin{align}\label{G12-1}
\Bigg| \frac{|\u|^{-3}}{H_\delta(\sin(\eps/2)) }  \Big(
 4 i(\u\cdot\zeta) H_\delta(1) &-  2(\u\cdot\zeta)^2 (C_\delta(1)- C_\delta(\sin(\eps/2))) \\
&- |\zeta^\perp|^2|\u|^2 (H_\delta(1) - C_\delta(1) +  C_\delta(\sin(\eps/2)) )\Big) \Big|  \nonumber\\ 
&\quad\qquad\qquad \  \  \le \frac{12\Gamma}{|H_\delta(\sin(\eps/2))|}\, \Big(\frac{|\zeta|}{|\u|^2 } + \frac{|\zeta|^2}{|\u| } \Big) \, . \nonumber
\end{align}
as defined in \eqref{HC-conditions}.

In order to control the remainder term $G_{b_\eps^\delta,3}$,  we write
%
 \begin{align}\label{Gn-22}
G_{b_\eps^\delta,3} =&  |\u|^{-3} \int_0^\pi \int_0^{2\pi} b_\eps^\delta (\cos \theta)  \sin\theta \left( \frac{-i e^{ic}}{6} \right) \\
& \times  \left( (\u \cdot \zeta)\sin^2(\theta/2) - |\u||\zeta^\perp| \sin(\theta/2) \cos(\theta/2) \sin\phi\right) ^3  d\phi d\theta \nonumber \\
=&  \frac{i |\u|^{-3}}{12 \pi H_\delta(\sin(\eps/2)) } \int_\eps^\pi \int_0^{2\pi} b^\delta (\cos \theta)  \sin\theta e^{ic} |\u|^3 |\zeta|^3 \sin^3(\theta/2)  \nonumber \\
& \times  \left( \cos\alpha \sin(\theta/2) - \sin\alpha \cos(\theta/2) \sin\phi\right) ^3  d\phi d\theta \nonumber 
 \end{align}
 where $\alpha$ is the angle between $\zeta$ and $\u$. 
 
 Thus, we can bound $G_{b_\eps^\delta,3}$ by
\begin{align}\label{Gn-23}
 |G_{b_\eps^\delta,3}| \le & \frac{|\zeta|^3}{12\pi | H_\delta(\sin(\eps/2)) |} \int_\eps^\pi \int_0^{2\pi} b^\delta(\cos\theta) \sin\theta \sin^3(\theta/2) \\
 &\times |\cos\alpha \sin(\theta/2) - \sin\alpha \cos(\theta/2) \sin\phi|^3 d\phi d\theta \nonumber \\
 \le &  \frac{|\zeta|^3}{6 |H_\delta(\eps)|} \int_\eps^\pi b^\delta(\cos\theta) \sin\theta \sin^3(\theta/2) \nonumber \\
 &\times \left[\sin^3(\theta/2) + 3\sin^2(\theta/2) \cos(\theta/2) + 3\sin(\theta/2)\cos^2(\theta/2) + \cos^3(\theta/2)\right] d\theta\nonumber\\
  \le &  \frac{4|\zeta|^3}{3 |H_\delta(\sin(\eps/2)) |} \int_\eps^\pi b^\delta(\cos\theta) \sin\theta \sin^3(\theta/2) d\theta. \nonumber
  \end{align}
  
Since, from \eqref{H-funct},  $H_\delta'(x)dx = \frac14 b^\delta(\cos\theta) \sin\theta \sin^2(\theta/2) d\theta$, we can explicitly calculate the error as a function of the
singular behavior of $b^\delta(\cos\theta)$.

\

 
 We first focus on the Rutherford-like cross section case ($\delta = 0$)  where,  by \eqref{HC-funct-2},  $H_0(x)=\log x$. Then
 \begin{align}\label{Gn-3.0}
 & H_0(\sin(\eps/2)) =  \log(\sin(\eps/2))\, , \ \ \textrm{and} \\ 
 & \int_\eps^\pi b^0 (\cos\theta)\sin\theta \sin^3(\theta/2) d\theta = 4\int_{\sin(\eps/2)}^1 \log x  dx = (1-\sin(\eps/2)) \le 4\,  . \nonumber
  \end{align}
  
 For the $0 < \delta < 2,\, \delta \neq 1$ case we have, by \eqref{H-cross-4},  $H_\delta(x)=-\frac{x^{-\delta}}{\delta}$. Then 
 \begin{align}\label{Gn-3.1}
& -H_\delta (\sin(\eps/2)) = \frac{\sin^{-\delta}(\eps/2)}{\delta}\, , \ \ \textrm{and}  \nonumber \\
& \int_\eps^\pi b^\delta(\cos\theta)\sin\theta \sin^3(\theta/2) d\theta =  -\frac{4}{\delta}\int_{\sin(\eps/2)}^1  x^{-\delta} dx= \\
&\qquad\qquad \qquad\qquad -\frac4{\delta(1-\delta)} (1 - \sin^{1-\delta} (\eps/2) ) \ \le\  \frac4{1-\delta} \, . \nonumber
 \end{align}
 
Thus,  from  \eqref{Gn-3.0} and \eqref{Gn-3.1},  we obtain the remainder term 
\begin{equation}\label{Gn-4.1}
G_{b_\eps^\delta,3}  \le O\left( \frac{|\zeta|^3}{H_\delta(\sin(\eps/2))} \right)\,  \ \ \textrm{for}\ \ 0\le\delta < 2,\, \delta\neq 1\,.
\end{equation}

The case  $\delta=1$  (so called $\eps$-linear case)  is special.  Indeed, by   \eqref{H-cross-4}, $H_1(x)=-x^{-1}$. Then 
 \begin{align}\label{Gn-3.2}
 & -H_1(\sin(\eps/2)) = \sin^{-1}(\eps/2)\, , \ \ \textrm{and}  \nonumber \\
 & \int_\eps^\pi b^1(\cos\theta)\sin\theta \sin^3(\theta/2) d\theta =  -4\int_{\sin(\eps/2)}^1  x^{-1} dx \\
 &\qquad= -4 \log (\sin(\eps/2)) \le 4 \left |\log (\sin(\eps/2))\right | \, , \nonumber
 \end{align}
and so 
corresponding  remainder term is of order
\begin{equation}\label{Gn-4.2}
G_{b_\eps^\delta,3}  \le O\left( 4 \left |\log (\sin(\eps/2))\right | \,   |\zeta|^3 \right)\, ,   \ \ \textrm{for}\ \ \delta=1.
\end{equation}

  
%

\noindent Gathering the estimates from \eqref{G12-1} and \eqref{Gn-23}, with \eqref{Gn-4.1} and \eqref{Gn-4.2}, 
 we obtain
\begin{align}\label{LandauConverge0}
|\tilde{G}&(\zeta,\u)| \le\ 
 \frac{12}{ \left |H_\delta(\sin(\eps/2)) \right | }\left( \Big(\frac{|\zeta|}{|\u|^2 } + \frac{|\zeta|^2}{|\u| }\Big) \Gamma  +  O(  {|\zeta|^3} ) \right) \, 1_{(0\le\delta<2, \delta \neq 1)} \\ 
&+12 \left |\sin(\eps/2) \right | \Big(\frac{|\zeta|}{|\u|^2 } + \frac{|\zeta|^2}{|\u| }\Big) \Gamma +  \left |\sin(\eps/2)\, \log (\sin(\eps/2))\right | 1_{\delta=1}  \,  O(  {|\zeta|^3} )
\, . \nonumber 
\end{align}


\medskip

\noindent So, for $0\le\delta<2,\, \delta \neq 0$  and $\eps<1$
\begin{align} \label{LandauConverge2}
 & |\widehat{Q_{b_\eps^\delta}}[f_\eps^\delta](\zeta) - \widehat{Q_L}[f_\eps^\delta](\zeta)| \ \le  \left(\frac{ 1_{\delta=0}}{\left | \log(\sin(\eps/2)) \right | }  + \frac{ 
 \left | \sin^\delta(\eps/2) \right | }{\delta(1-\delta)} 1_{(0<\delta<2,\, \delta\neq 1)}\right) \\
&\qquad\qquad \qquad\qquad   \Bigg| \int_{\R^3} \mathcal{F}\{f_\eps^\delta(v)f_\eps^\delta(v-u)\}(\zeta) 
\, O\left( \Big(\frac{|\zeta|}{|\u|^2} + \frac{|\zeta|^2}{|\u| }\Big) + |\zeta|^3 \right) du \Bigg|\, , \nonumber 
 \end{align}
and for $\delta=1$ and $\eps<1$
\begin{align} \label{LandauConverge3}
 & |\widehat{Q_{b_\eps^\delta}}[f_\eps^\delta](\zeta) - \widehat{Q_L}[f_\eps^\delta](\zeta)| \ \le \  \left | \sin(\eps/2)\, \log (\sin(\eps/2)) \right |  \\
&\qquad  \Bigg| \int_{\R^3} \mathcal{F}\{f_\eps^\delta(v)f_\eps^\delta(v-u)\}(\zeta) \, O\left( \Big(\frac{|\zeta|}{|\u|^2} + \frac{|\zeta|^2}{|\u| }\Big) +
 |\zeta|^3 \right) du \Bigg| .\nonumber
 \end{align}
 \ 

\noindent Finally,  using assumption \eqref{fAssump} on the Fourier transform of $f_\eps^\delta(\v)f_\eps(\v-\u)$ yields the sufficient decay for integrability at infinity, 
making the above integral finite. Thus,  
\begin{align}\label{LandauConverge-2}
 &\Big|\widehat{Q_L}[f_\eps^\delta](\zeta) - \widehat{Q_{b_{\eps}}}[f_\eps^\delta](\zeta) \Big| \  \le   \nonumber  \\
 &\ \\
	&\ \  \ O\left( \frac{ 1_{\{\delta=0\}}} {| \log(\sin(\eps/2))| } +
	  \frac{ \left | \sin^{\delta}(\frac{\eps}2) \right | }{\delta(1-\delta)} 1_{\{0<\delta<2,\,\delta \neq 1\}} 
	+ \left | \sin(\eps/2)\, \log (\sin(\eps/2)) \right | 1_{\{\delta=1\}} \right)  \nonumber \\
		&\    \nonumber \\
&=  
O\left(\frac{  \left | 1 +  (| \log(\sin(\eps/2))|  - 1)  \, 1_{\{\delta=1\}} \right |} { \left | H_\delta(\sin(\eps/2)) \right |}
\right)\ \rightarrow 0 \, ,\nonumber
\end{align}
uniformly in $\zeta\in \R^3$ for each fixed $\delta$, and so \eqref{FTAssump}  holds.
\qed
\end{proof}

\bigskip

The result from this theorem  illustrates that \emph{the convergence to the Landau collision operator is highly dependent on the model chosen for $b_\eps^\delta(\cos \theta)$.}

 The case of Rutherford scattering corresponds to the choice of $H$ from \eqref{HC-funct-2}, which results in logarithmically slow convergence to the Landau operator \eqref{FTAssump}. As was noted in Remark 1, the rescaling of the cross section that was done to construct $b_\eps$ ensured that the leading order term in the expansion is $O(1)$. Without this rescaling, the leading order term is the Landau operator scaled by $\log(\sin(\eps/2))$, i.e., the Coulomb logarithm that appears in the original derivation by Landau, and in this case the next order term in the expansion would give an $O(1)$ error.
 

\medskip

\begin{remark}[Remark 4:]\label{rem4} This rescaling implies that the relevant time scale for grazing collision effects with the Coulombic potential is much longer than the mean free time between collisions typically used in non-dimensionalizing the Boltzmann equation. When solving a space inhomogeneous problem, this presents a time scale separation between the grazing collisional terms and the Vlasov terms on the left hand side of a space inhomogenous problem. Solving Coulombic interaction problems using the Boltzmann equation with finite $\eps$ therefore captures both the faster scale strong collisions as well as the weak collisions. Furthermore, using Boltzmann gives a faster decay rate to equilibrium than the Landau operator alone owing to inclusion of the error term controlled by $\log \sin(\eps/2)$  (as shown in Figure~\ref{LvsCEnt}). However, it should be noted that while Boltzmann collision operator with small $\eps$ includes the effects of both strong collisions and grazing collisions, it is missing another key effect that comes into play at small $\eps$, namely, the influence of collective effects at large interaction distance that is modeled by the Lenard-Balescu collision operator \cite{FriBook}.
\end{remark}

%

\medskip


\bigskip



\section{The Conservative Numerical Method} \label{sec:numerics_setup}

\subsection{Velocity space discretization} \label{sec:discretization}

In order to compute the Boltzmann equation we must work on a bounded velocity space, rather than all of $\R^3$. However typical distributions are supported on the entire domain, for example the Maxwellian equilibrium distribution. Even if one begins with a compactly supported initial distribution, each evaluation of the collision operator spreads the support by a factor of $\sqrt{2}$, thus we must use a working definition of an {\em effective support} of size  $R$ for the distribution function. Bobylev and Rjasanow \cite{BobRja99} suggested using the temperature of the distribution function, which typically decreases as $\text{exp}(-|v|^2 / 2T)$ for large $|v|$, and used a rough estimate of $R \approx 2\sqrt{2}T$ to determine the cutoff. Thus, we assume that the distribution function is negligible outside of a ball 
\begin{equation} \label{Ball_domain_v}
B_{R_x}(\textbf{V}(\x)) = \{\v \in \R^3 : |\v - \textbf{V(\x)}| \le R_x \},
\end{equation}
where $\textbf{V}(\x)$ is the local flow velocity which depends in the spatial variable $\x$. For ease of notation in the following we will work with a ball centered at $0$ and choose a length $R$ large enough that $B_{R_x}(\textbf{V}(\x)) \subset B_R(0)$ for all $\x$.

With this assumed support for the distribution $f$, the integrals in \eqref{Cont_spectral} will only be nonzero for $\u \in B_{2R}(0)$. Therefore, we set $L=2R$ and define the cube
\begin{equation} \label{Cube_domain_v}
C_L = \{ \v \in \R^3 : |v_j| \le L,\,\, j = 1,\dots,d\}
\end{equation}
to be the domain of computation. With this domain the computation of the weight function integral \eqref{Ghat_Aniso} is cut off at $r=L$.

Let $N \in \mathbb{N}$ be the number of points in velocity space in each dimension. Then we establish a uniform velocity mesh with $\Delta v = \frac{2L}{N-1}$ and due to the formulation of the discrete Fourier transform the corresponding uniform Fourier space mesh size is given by $\Delta \zeta = \frac{(N-1)\pi}{NL}$. 
%


\subsection{Collision step discretization} \label{sec:collision_disc}
%

To simplify notation we will use one index to denote multidimensional sums with respect to an index vector $\textbf{m}$
\begin{equation*}
\sum_{\textbf{m}=0}^{N-1} = \sum_{m_1,\dots,m_d = 0}^{N-1}.
\end{equation*}

To compute $\widehat{Q}(\zeta_\k)$, we first compute the Fourier transform integral giving $\hat{f}(\zeta_k)$ via the FFT. The weighted convolution integral is approximated using the trapezoidal rule
\begin{equation}
\widehat{Q}(\zeta_\k)= \sum_{\textbf{m} = 0}^{N-1} \widehat{G}(\xi_{\textbf{m}},\zeta_\k) \hat{f}(\xi_{\textbf{m}}) \hat{f}(\zeta_\k - \xi_{\textbf{m}}) \omega_{\textbf{m}},
\end{equation}
where $\omega_\textbf{m}$ is the quadrature weight and we set $\hat{f}(\zeta_\k - \xi_{\textbf{m}}) = 0$ if $(\zeta_\k - \xi_{\textbf{m}})$ is outside of the domain of integration. We then use the inverse FFT on $\widehat{Q}$ to calculate the integral returning the result to velocity space. 

Note that in this formulation the distribution function is not periodized, as is done in the collocation approach of Pareschi and Russo \cite{ParRus00}. This is reflected in the omission of Fourier terms outside of the Fourier domain. All integrals are computed directly only using the FFT as a tool for faster computation and the convolution integral is accurate to at least the order of the quadrature. The calculations below use the trapezoid rule, but in principle Simpson's rule or some other uniform grid quadrature can be used. However, it is known that the trapezoid rule is spectrally accurate for periodic functions on periodic domains (which is the basis of spectral accuracy for the FFT), and the same arguments can apply to functions with sufficient decay at the integration boundaries \cite{Atkinson}. These accuracy considerations will be investigated in future work. The overall cost of this step is $O(N^{6})$. 

\subsection{Discrete conservation enforcement} \label{sec:conservation}
This implementation of the collision mechanism does not conserve all of the quantities of the collision operator. To correct this, we formulate these conservation properties as a Lagrange multiplier problem. Depending on the type of collisions we can change this constraint set (for example, inelastic collisions do not preserve energy), but we will focus on the case of elastic collisions, which preserve mass, momentum, and energy. 

Let $M = N^d$ be the total number of grid points, let $\tilde{\textbf{Q}} = (\Qt_1, \dots, \Qt_M) ^T$ be the result of the spectral formulation from the previous section, written in vector form, and let $\omega_j$ be the quadrature weights over the domain in this ordering. Define the integration matrix
\begin{equation*}
\textbf{C}_{5\times M} = \left(\begin{array}{c} \omega_j \\ v_j^i \omega_j \\ |\v_j|^2 \omega_j \end{array} \right),
\end{equation*}
where $v^i,\, i=1,2,3$ refers to the $i$th component of the velocity vector. Using this notation, the conservation method can be written as a constrained optimization problem. 

\begin{equation} 
\text{Find } \textbf{Q} = (Q_1,\dots,Q_M)^T \text{ that minimizes } \frac12 \|\tilde{\textbf{Q}} - \textbf{Q}\|_2^2 \text{ such that } \textbf{C} \textbf{Q} = \textbf{0}
\end{equation}
Formulating this as a Lagrange multiplier problem, we define 
\begin{equation}
L(\textbf{Q},\gamma) = \sum_{j=1}^M (\Qt_j - Q_j)^2 - \gamma^T\textbf{C}\textbf{Q}
\end{equation}
The solution is given by
\begin{align}
\textbf{Q} &= \tilde{\textbf{Q}} + \textbf{C}(\textbf{C}\textbf{C}^T)^{-1} \textbf{C} \tilde{\textbf{Q}} \nonumber \\
&:= \textbf{P}_N \tilde{\textbf{Q}}
\end{align}

Overall the collision step in semi-discrete form is given by
\begin{equation}
\frac{\d \textbf{f}}{\d t} = \textbf{P}_N \tilde{\textbf{Q}}
\end{equation}

The overall cost of the conservation portion of the algorithm is a $O(N^d)$ matrix-vector multiply, significantly less than the computation of the weighted convolution.

\subsection{Computing $\Ghat$ for singular scattering cross sections}

Numerically calculating the weights $\hat{G}$ to high accuracy can be difficult for singular scattering cross sections, due to the precise nature of the cancellation at the left endpoint of the integral. The $\theta$ integral in \eqref{Ghat_Aniso} can be simplified as
\begin{equation}
\int_\eps^\pi b_\eps (\cos \theta) \sin \theta (\cos(c_1 (1-\cos\theta) - c_3)  J_0(c_2 \sin \theta) - \cos(c_3)) d\theta,
\end{equation}
where $c_1, c_2, c_3$ depend on the current values of $\phi, r, \zeta, \xi$ following from the full formulation of $\Ghat$. When $\eps << 1$ the bulk of the integration mass occurs near the left endpoint of the $\theta$ interval, however this presents a challenge for a numerical quadrature package to compute. {For $\theta << 1$ there is a subtraction of two nearly equal numbers (the two cosine terms), which causes floating point errors. To alleviate this, we split the integration interval into two pieces, and use the first term of the Taylor expansion of the troublesome part of the integrand for $\theta << 1$ and obtain
\begin{align*}
&\left(-\frac{c_2^2}{4}\cos(c_3) + \frac{c_1}{2}\sin(c_3)\right) \int_\eps^{\sqrt{\eps}} \theta^2 b_\eps(\cos \theta) \sin \theta d\theta\ \nonumber \\
+&\int_{\sqrt{\eps}}^\pi b_\eps(\cos \theta) \sin \theta (\cos(c_1 (1-\cos\theta) - c_3)  J_0(c_2 \sin \theta) - \cos(c_3)) d\theta.
\end{align*}
}
These integrals are computed using the GNU Scientific Laboratory integration routines \cite{gsl}. We use \verb+cquad+ to compute the first $\theta$ integral, which appears to be most stable choice for this near-singular integrand. The adaptive Gauss-Konrod quadrature \verb+qag+ is used for all other integrals used in computing the weights $\Ghat$. This calculation is very expensive, however speedup of this high-dimensional calculation is done using OpenMP and MPI  on a cluster, as each weight can be calculated independently in parallel \cite{GamHaaRGD, HaackHPC}. Once the weights are known, the computational cost is the same as the original spectral method \cite{GamTha09, GamTha10}.
\section{Numerical results}
To illustrate that this method captures the correct behavior for grazing collisions, we take $\lambda = -3$ and set $\eps = 10^{-4}$. Similar to what was done in \cite{BuetCord, ParRusTos00} based on the original work of Rosenbluth et al. \cite{Rosenbluth57} for the Landau equation, we set the axially symmetric initial condition
\[ f(\v,0) = 0.01  \textrm{exp}\left(-10 \left(\frac{|\v| - 0.3}{0.3}\right)^2\right). \]

\begin{figure}[ht] 
\centering
\includegraphics[width=.95\linewidth]{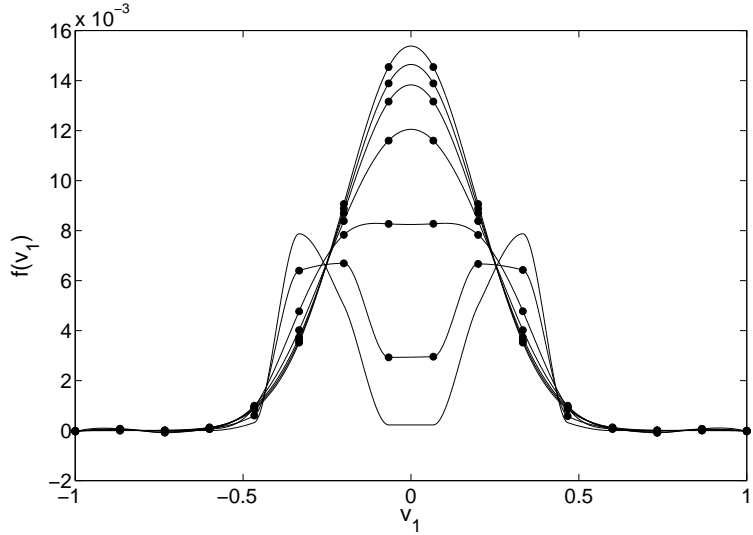}
\caption{Slice of the distribution marginal function at times $t=0, 9, 36, 81, 144, 225, 900$. Solid lines: Hermite spline reconstruction of Landau equation solution. Solid circles: Boltzmann solution with artificial cross section \eqref{beps} . $\eps = 10^{-4}, N=16$.}
\label{GrazingSlice}
\end{figure}

We begin by using the $\eps$-linear cross section \eqref{beps} instead of the Rutherford cross section, as the logarithmic error term is too large to make any meaningful convergence observations, even if we take $\eps$ to machine precision. Note that from the analysis in section 3, we expect the error to be $O\left(\frac{\log\sin((\eps/2))}{H_1(\sin(\eps/2))}\right) = O(\sin(\eps/2)\log(\sin(\eps/2))) \sim O(\eps \log \eps)$, essentially $O(\eps)$ for the values of $\eps$ we use. We take a domain size of $L = 1$, glancing parameter $\eps = 10^{-4}$, $N=16$, and compute to time $t = 900$ with a timestep of $0.01$. The results are shown in Figure \ref{GrazingSlice}. Note that our symmetric grid is not aligned with $v_1=0$, so it is slightly offset from the figures from the earlier works.

\begin{figure}[htb] 
 \begin{center}
\subfigure[Early times]{\includegraphics[width=4.5in,angle=0]{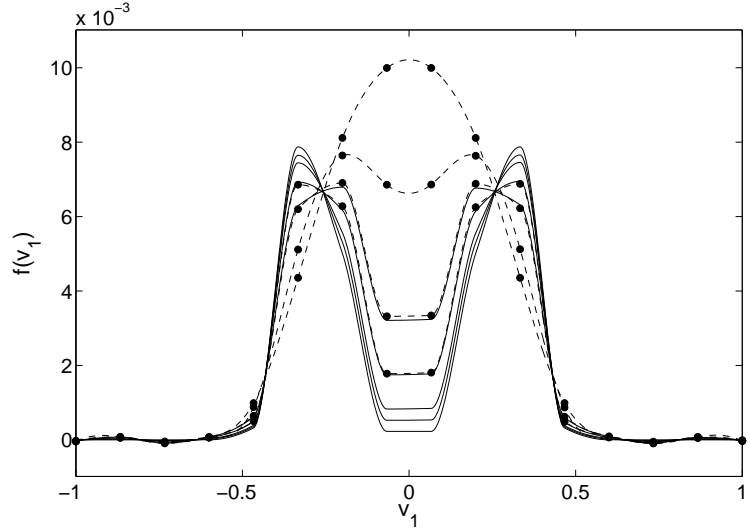}} 
\subfigure[Long times]{\includegraphics[width=4.5in,angle=0]{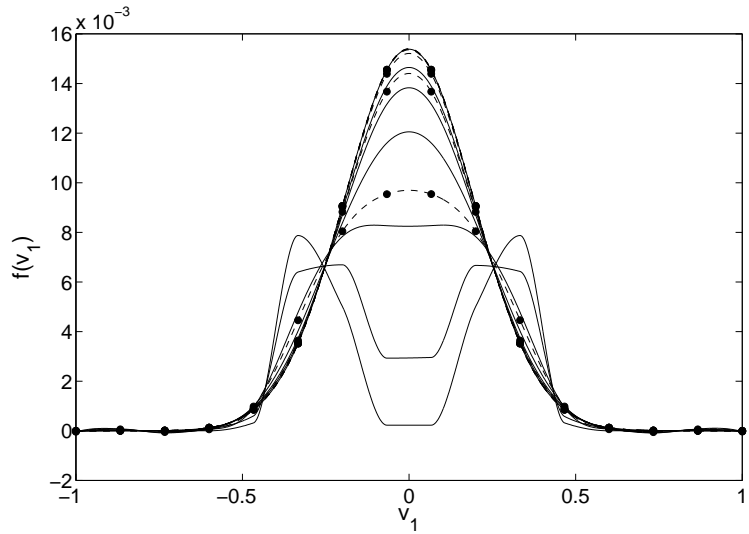}}
\end{center}
\caption{{Comparisons of solutions to Boltzmann using Rutherford cross section \eqref{beps_ruth} and to Landau equations.
{\sl (a)}  Slice of the distribution marginal function at early times $t=0, 1, 2, 5, 10$. {\sl (b)} Slice of the distribution marginal function at times $t=0,  9, 36, 81, 144, 225, 900$.  Solid lines: spline reconstruction of Landau equation solution. Dashed lines with solid circles: spline reconstruction of Boltzmann equation. Spline reconstruction uses Hermite polynomials for times below $t=10$ to avoid a reconstruction the generate negative values in the marginal tails $\eps = 10^{-4}, N=16$.} }
\label{LvsCT10}
\end{figure}

To verify the linear convergence rate for the artificial cross section, we take a single timestep of the example above for $\eps=10^{-1},10^{-2},10^{-3},10^{-4}$ with the artificial cross section \eqref{beps}. We compare these values with the result of a single step of the Landau equation, computed using the spectral method with the convolution weights derived for the Landau equation \eqref{GFPL}. We represent the error by examining the difference in the values in the central slice of the solution, which are the same values plotted in Figures \ref{GrazingSlice} and \ref{LvsCT10}. In Table \ref{times} we present the average error between the Landau and Boltzmann solution in this subset. As expected, the convergence is linear.

%
%

\begin{table}[!htbp] 
\begin{tabularx}{\textwidth}{XXXX}
\hline
$\eps$ & average $|Q_L - Q_\eps|$ & ratio\\
\hline
$10^{-1}$ & $1.35 \times 10^{-4}$ & \\
$10^{-2}$ & $1.47 \times 10^{-5}$ & 8.98\\
$10^{-3}$ & $1.51 \times 10^{-6}$ & 9.68\\
$10^{-4}$ & $1.52 \times 10^{-7}$ & 9.89\\
\hline
\end{tabularx}
\caption{Error between Boltzmann collision operator with grazing collisions and Landau collision operator}
\centering
\label{times}
\end{table}

\begin{figure}[ht]
\centering
\includegraphics[width=\linewidth]{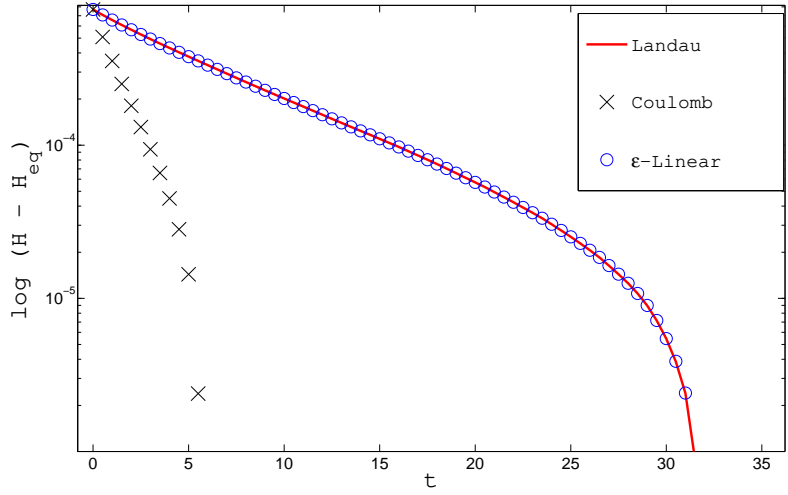}
\caption{Convergence of entropy to equilibrium: Log of entropy decay for Boltzmann solution with the Rutherford cross section \eqref{beps_ruth}  with  crosses, Boltzmann solution with the {$\eps$-linear} cross section \eqref{beps}   with  circles, and  Landau solution with solid curve. $N=16,\,\eps=10^{-4}$. When calculating the entropy $H$, we exclude grid points where the distribution is negative.}
\label{LvsCEnt}
\end{figure}

In Figure \ref{LvsCT10} we plot the evolution of the Boltzmann equation using the Rutherford cross section \eqref{beps_ruth} and compare it to the numerical solution of the Landau equation. We again take $L=1,\, \eps=10^{-4}$, and $N=16$. This figure illustrates the large error between the two models for this cross section, as well as the different convergence rates to equlibrium. Indeed, we can see this more explicitly in Figure \ref{LvsCEnt}, where we can see the solution of the Boltzmann equation converges to equilibrium at a much faster rate than the Landau equation. We also ran tests for much smaller values of $\eps$ with this cross section (e.g. $\eps = 10^{-12}$), however due to the slow rate of convergence of the error term this did not result in a significant difference in the solution.

Here we remark that the recent work of Bobylev and Potapenko \cite{BobPot} proves that the order of approximation between the Boltzmann and Landau operators is no worse than $\sqrt{\eps}$, which would seem to contradict our logarithmic convergence result. However, the effective cross section used in their work does not satisfy the assumptions on the scattering cross section in the theory above, so there is no contradiction.

\begin{figure}[ht] 
\centering
\includegraphics[width=3in,angle=0]{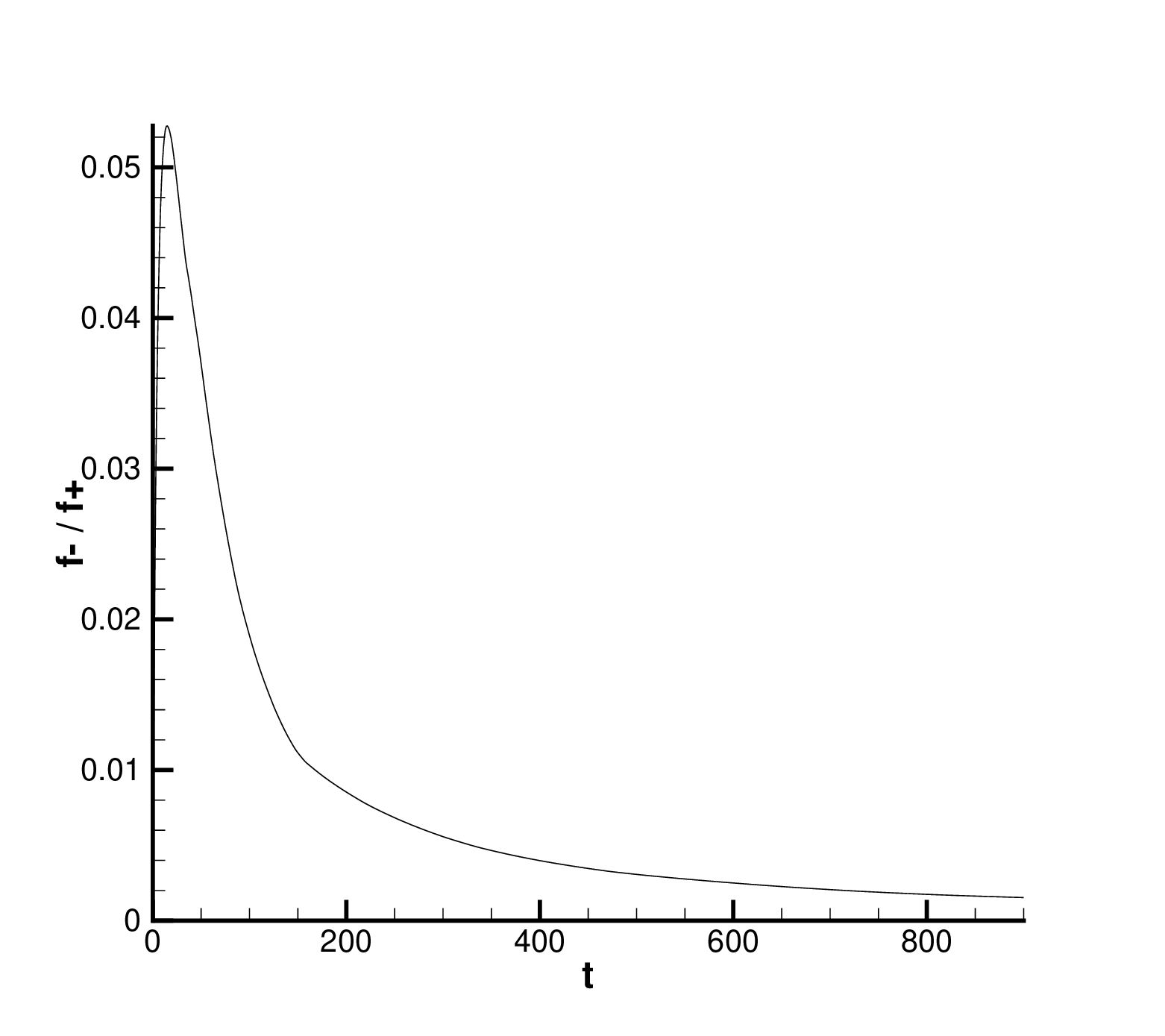}
\caption{Ratio of energy in negative grid points to energy in positive grid points. $\eps = 10^{-4}, N=16$}
\label{GrazingEnRat}
\end{figure}

Due to the spectral formulation some of the grid values may be negative. Recent work by Alonso, Gamba, and Tharkabhushanam \cite{AlonGamTha} has shown that the scheme maintains its spectral properties, converges and obtain error estimates provided that the `energy' of the negative grid points remains small compared to the energy of the rest of the computed distribution. In Figure \ref{GrazingEnRat} we plot the percentage 
\[ \frac{\int_D |f_{-}| |v|^2 dv}{\int_D f_{+} |v|^2 dv}, \]
where $f_{-}$ is the grid cells where the discrete distribution function is negative and vice versa for $f_{+}$. Future work will explore including positivity constraints to the conservation routines.


\section{Conclusions and future work}
We have derived the spectral formulation for the more general case of anisotropic collisional models for the Boltzmann equation. We also showed that the spectral method for the Boltzmann equation is consistent with the limiting Landau equation under suitable assumptions on the scattering cross section, and that using the grazing collision Boltzmann equation can give very different results in convergence to the Landau equation depending on the cross section chosen. We also note that the grazing collisions limit was chosen only because it was a case in which we have an explicit model to compare our anisotropic Boltzmann solver with. If one wants to solve the Landau equation, one can simply use the weights \eqref{GFPL} derived from the weak form of the Landau equation within the spectral framework without worrying about the angular dependence. One other important thing to note is that this method may be a good candidate for collisional models where the collision mechanism is unknown and only experimentally determined, and future work will attempt to simulate the Boltzmann equations with these cross sections. In addition, as the Landau equation is used to model collisions of charged particles in plasma we will seek to add field effects to the space inhomogeneous Boltzmann equation, resulting in the Boltzmann-Poisson or Boltzmann-Maxwell systems. The inhomogeneous method uses operator splitting between the collision and the transport terms, so in principle one can use an already developed Vlasov solver for the spatial terms in the equation. 

\section*{Acknowledgments}
The authors thank the referees for their helpful comments. The authors would also like to thank Phil Morrison for many helpful discussions about plasma physics. This work has been supported by the NSF under grants DMS-0636586,  DMS-1109625, and NSF RNMS (KI-Net) grant \#11-07465.  





\bibliographystyle{siam}   

\bibliography{Boltz}


\end{document}